\documentclass[12pt]{article}
\usepackage{amssymb,amsmath}
\usepackage{cases}
\usepackage{amsfonts}
\usepackage{color,xcolor}
\usepackage[left=2.0cm,right=2.0cm,top=2.0cm,bottom=2.0cm]{geometry}
\usepackage[colorlinks,citecolor=blue,urlcolor=blue]{hyperref}

\newtheorem{theorem}{Theorem}[section]

\newtheorem{lemma}{Lemma}[section]

\newtheorem{remark}{Remark}[section]

\newcommand{\bal}{\begin{align}}
\newcommand{\bbal}{\begin{align*}}
\newcommand{\beq}{\begin{equation}}
\newcommand{\eeq}{\end{equation}}
\newcommand{\bca}{\begin{cases}}
\newcommand{\eca}{\end{cases}}
\def\div{\mathord{{\rm div}\ }}
\newcommand{\pa}{\partial}
\newcommand{\fr}{\frac}
\newcommand{\na}{\nabla}
\newcommand{\De}{\Delta}

\newcommand{\cd}{\cdot}
\newcommand{\lan}{\langle}
\newcommand{\ran}{\rangle}
\newcommand{\ep}{\varepsilon}
\newcommand{\dd}{\ \mathrm{d}}

\newcommand{\R}{\mathbb{R}}

\newcommand{\bi}{\Big}

\begin{document}
\title{A class of large solution of 2D tropical climate model without thermal diffusion}

\author{Jinlu Li$^{1}$\footnote{E-mail: lijinlu@gnnu.cn}, Yanghai Yu$^{2}$\footnote{E-mail: yuyanghai214@sina.com(Corresponding author)} \\
\small $^1$\it School of Mathematics and Computer Sciences, Gannan Normal University, Ganzhou 341000, China\\
\small $^2$\it School of Mathematics and Statistics, Anhui Normal University, Wuhu, Anhui, 241002, China}

\date{}

\maketitle\noindent{\hrulefill}

{\bf Abstract:} In this paper, we consider the Cauchy problem of 2D tropical climate model without thermal diffusion and construct global smooth solutions by choosing a class of special initial data whose $L^{\infty}$ norm can be arbitrarily large.


\vskip0mm\noindent{\hrulefill}

\section{Introduction}\label{sec1}
This paper focuses on the following 2D tropical climate model (TCM) given by
\begin{eqnarray}\label{2d-tcm}
        \left\{\begin{array}{ll}
          \partial_tu+u\cd\na u+\mu u+\na p+\div(v\otimes v)=0,& x\in \R^2,t>0,\\
          \partial_tv+u\cd\na v+v\cd\na u+\na \theta-\nu\Delta v=0,& x\in \R^2,t>0,\\
          \partial_t\theta+u\cd\na \theta+ \div v=0,& x\in \R^2,t\geq0,\\
          \div u=0,& x\in \R^2,t\geq0,\\
                   (u,v,\theta)|_{t=0}=(u_0,v_0,\theta_0),& x\in \R^2.\end{array}\right.
        \end{eqnarray}
here $\mu,\nu$ are non-negative parameters, $u=u(t,x)$ and $v=v(t,x)$ stand for the barotropic mode and the first baroclinic mode of the vector velocity, respectively. $p=p(t,x)$ and $\theta=\theta(t,x)$ denote the scalar pressure and scalar temperature, respectively.

In the completely inviscid case, namely, $\mu=\nu=0$, (\ref{2d-tcm}) was first derived by Frierson--Majda--Pauluis \cite{Frierson 2004} by performing a Galerkin truncation to the hydrostatic Boussinesq equations, of which the first baroclinic mode had been originally used in some studies of tropical atmosphere. More relevant background on the tropical climate model can be found in \cite{Gill 1980,Matsuno 1966,Biello 2003} and the references therein. The mathematical studies on (\ref{2d-tcm}) have attracted considerable attention recently from various authors and have motivated a large number of research papers concerning the local well-posedness \cite{Ma 2016}, global small solutions, global regularity and so on. Li--Titi \cite{Li 2016} introduced a new quantity to bypass the obstacle caused by the absence of thermal diffusion and proved the global well-posedness of the 2D TCM (\ref{2d-tcm}) with two Laplacian terms $\Delta u$ and $\Delta v$. From the mathematical point of view, adding some fractional dissipation terms does bring the regularity for the system (\ref{2d-tcm}) and also significantly change the model properties and physics. Concerning the global regularity of TCM with fractional dissipation terms, we refer to \cite{Dong1 2018,Dong2 2018,Dong3 2018,Ye 2017,Zhu 2018} and the references therein. However, adding only the damping term $u$ does not appear to be enough to establish the global well-posedness of TCM (\ref{2d-tcm}). Wan \cite{Wan 2016} proved the global strong solution to TCM (\ref{2d-tcm}) with additional damping term $v$ under suitable smallness assumptions on the initial data. Later, Ma--Wan \cite{Ma 2018} removed the damping term $v$ and obtained the global strong solution to TCM (\ref{2d-tcm}) with different smallness assumptions on the initial data. Naturally, we ask that whether or not there exists a global solution to TCM (\ref{2d-tcm}) with some classes of large initial value. Motivated by the ideas that used in \cite{Liyz,Liy,Lj 2019}, we give a positive answer in this paper.

For the sake of implicity, we will limit ourselves to $\mu=\nu=1$  since their values do not play any role in our analysis and introduce the following notations
\bbal
A(t)\triangleq||(w,c,\theta)(t)||^2_{H^3}, \quad B(t)\triangleq||(w,\na c)(t)||^2_{H^3}+||\na\theta(t)||^2_{H^2}\quad\mbox{and}\quad A_0\triangleq||w_0,c_0,\theta_0||^2_{H^3}.
\end{align*}
Our main result is stated as follows.
\begin{theorem}\label{the1.1} Assume that the initial data fulfills $\div w_0=0$ and
$$u_0=U_0+w_0\quad \mbox{and}\quad v_0=V_0+c_0$$
where
\begin{eqnarray*}
&U_0=
\begin{pmatrix}
\pa_2a_0 \\ -\pa_1a_0
\end{pmatrix}
\quad\mbox{and}\quad
V_0=
\begin{pmatrix}
m_0 \\ m_0
\end{pmatrix}
\end{eqnarray*}
with
\begin{eqnarray}\label{con-Equ1.2}
\mathrm{supp} \ \hat{a}_0(\xi),\mathrm{supp}\ \hat{m}_0(\xi)\subset\mathcal{C}\triangleq\Big\{\xi\in\R^2: \ |\xi_1+\xi_2|\leq \ep, \ 1\leq |\xi|\leq 2\Big\},\quad \ep>0 .
\end{eqnarray}
There exists a sufficiently small positive constant $\varepsilon_0$, and a universal constant $C$ such that if
\begin{align}\label{condition}
&\Big(A_0+\ep\big(||a_0||_{L^2}||\hat{a}_0||_{L^1}+||m_0||_{L^2}(1+||\hat{m}_0||_{L^1})\big)\Big)\exp\Big\{C\ep\big(||a_0||_{L^2}||\hat{a}_0||_{L^1}+||m_0||_{L^2}(1+||\hat{m}_0||_{L^1})\big)
\nonumber\\&\quad~~~~~~~~~~~~~~~~~ +C\big(||\hat{a}_0,\hat{m}_0||_{L^1}+||\hat{a}_0,\hat{m}_0||^2_{L^1}
\big)\Big\}\leq \varepsilon_0,
\end{align}
then the system \eqref{2d-tcm} has a unique global solution.
\end{theorem}
\begin{remark}\label{rem1.1} By choosing a class of special initial data, we can show that general initial value whose $L^{\infty}$ norms can be arbitrarily large generates a global unique solution to the system \eqref{2d-tcm}. Here, we construct an example to verify this.

Let $(w_0,c_0,\theta_0)=(0,0,0)$ and $a_0=m_0=\frac1\ep\big(\log\log\frac1\ep\big)^\frac12 \chi$, where the smooth function $\chi$ satisfying
\begin{align*}
\mathrm{supp} \hat{\chi}\in \mathcal{{C}},\quad \hat{\chi}(\xi)\in[0,1]\quad\mbox{and} \quad \hat{\chi}(\xi)=1 \quad\mbox{for} \quad \xi\in\mathcal{\widetilde{C}},
\end{align*}
where
\begin{align*}
\mathcal{\widetilde{C}}\triangleq\Big\{\xi \in\R^2:  \ |\xi_1+\xi_2|\leq \frac\ep2,\ \frac43\leq |\xi|\leq \frac53\Big\} .
\end{align*}
Then, direct calculations show that the left side of \eqref{condition} becomes
\begin{align*}
C\ep^\frac12\Big(\log\log \frac1\ep\Big)\exp\Big\{C\log\log \frac1\ep\Big\}.
\end{align*}
In fact, one has
\begin{align*}
||\hat{a}_0||_{L^1}\approx \Big(\log\log\frac1\ep\Big)^\frac12\quad\mbox{and}\quad||{a}_0||_{L^2}\approx \ep^{-\fr12}\Big(\log\log\frac1\ep\Big)^\frac12.
\end{align*}
Therefore, choosing $\ep$ small enough, we deduce that the system \eqref{2d-tcm} has a global solution.

Moreover, we also have
\begin{align*}
||u_0||_{L^\infty}\gtrsim||\De a_0||_{L^\infty}\gtrsim||\hat{a}_0||_{L^1}\gtrsim \bi(\log\log\frac1\ep\bi)^\frac12
\end{align*}
and
\begin{align*}
||v_0||_{L^\infty}\approx||m_0||_{L^\infty}\gtrsim||\hat{m}_0||_{L^1}\gtrsim \bi(\log\log\frac1\ep\bi)^\frac12.
\end{align*}
\end{remark}

\section{Renormalized System}\label{sec2}
Let $(a,m)=(e^{-t}a_0,e^{t\Delta}m_0)$ be the solutions of the following system
\begin{eqnarray}\label{l1}
        \left\{\begin{array}{ll}
          \pa_ta+a=0,\\
          \pa_tm-\De m=0,\\
          (a,m)|_{t=0}=(a_0,m_{0}).\end{array}\right.
        \end{eqnarray}
Setting
\bbal
&U=
\begin{pmatrix}
\pa_2a \\ -\pa_1a
\end{pmatrix}
 \quad\mbox{and}\quad
V=
\begin{pmatrix}
m \\ m
\end{pmatrix},
\end{align*}
we can deduce from \eqref{l1} that
\begin{eqnarray}\label{l2}
        \left\{\begin{array}{ll}
          \pa_t U+ U=0,\\
          \pa_t V-\De V=0,\\
          \div U=0,\\
          (U,V)|_{t=0}=(U_0,V_{0}).\end{array}\right.
        \end{eqnarray}
Denoting $w=u-U$ and $c=v-V$, the system \eqref{2d-tcm} can be written as follows
\begin{eqnarray}\label{c-tcm}
        \left\{\begin{array}{ll}
\partial_tw+w\cd\na w+U\cd\na w+c\cd\na c+c\div c+w+\na p=f\\
\qquad \qquad -c\cd\na V-c\div V-V\cd\na c-V\div c-w\cd\na U,\\
\partial_tc+w\cd\na c+U\cd\na c-\De c+\na\theta=g-c\cd\na w-V\cd\na w-c\cd\na U-w\cd\na V,\\
\partial_t\theta+w\cd\na \theta+U\cd\na \theta+ \div c=h,\\
\div w=0,\\
(w,c,\theta)|_{t=0}=(w_0,c_0,\theta_0).\end{array}\right.
\end{eqnarray}
where
\bbal
&f=-U\cd\na U-V\cd\na V-V\div V,\qquad  g=-U\cd\na V-V\cd\na U, \qquad h=-\div V.
\end{align*}
\section{Useful Tools}\label{sec3}
\setcounter{equation}{0}
Firstly, we introduce some notations and conventions which will be used throughout this paper.
\begin{itemize}
 \item We will use the notation $||f_1,\cdots,f_n||_{X}\triangleq||f_1||_{X}+\cdots+||f_n||_{X}$ for some Banach space $X$.
  \item Let $\alpha=(\alpha_1,\alpha_2)\in \mathbb{N}^2$ be a multi-index and $D^{\alpha}=\pa^{|\alpha|}/\pa^{\alpha_1}_{x_1}\pa^{\alpha_2}_{x_2}$ with $|\alpha|=\alpha_1+\alpha_2$.
  \item $\lan f,g\ran$ denotes the inner product in $L^2(\R^2)$, namely, $\lan f,g\ran\triangleq\int_{\R^2}fg\dd x$.
   \item The Fourier transform of $f$ with respect to the space variable is given by
$$\hat{f}(\xi)\triangleq\int_{\R^2}e^{-ix\cd\xi}f(x)dx.$$
  \item For $s\in \mathbb{N}$, the norm of the integer order Sobolev space $H^s(\R^2)$ and $W^{s,\infty}(\R^2)$ are defined by
  $$||f||_{H^s(\R^2)}\triangleq\bi(\sum_{|\alpha|\leq s}||D^{\alpha}f||^2_{L^2(\R^2)}\bi)^{\fr12}=||f||_{\dot{H}^s(\R^2)}+||f||_{L^2(\R^2)}$$
  and
   $$||f||_{W^{s,\infty}(\R^2)}\triangleq\sum_{0\leq i\leq s}||\nabla^{i}f||_{L^{\infty}(\R^2)}.$$
\end{itemize}

Next, we present some estimates which will be used in the proof of Theorem \ref{the1.1}.
\begin{lemma}\label{le2} \cite{Majda 2001} (Commutator estimates)
There hold that
\bal
&\sum_{0<|\alpha|\leq 3}||[D^{\alpha},\mathbf{f}]\mathbf{g}||_{L^2}\leq C(||\mathbf{g}||_{{H}^{2}}||\na \mathbf{f}||_{L^\infty}+||\mathbf{g}||_{L^\infty}||\mathbf{f}||_{{H}^3}),\label{jhz1}\\
&\sum_{0<|\alpha|\leq 3}||[D^{\alpha},\mathbf{f}]\mathbf{g}||_{L^2}\leq C(||\na \mathbf{f}||_{L^\infty}+||\na^3 \mathbf{f}||_{L^\infty})||\mathbf{g}||_{H^2}.\label{jhz2}
\end{align}
\end{lemma}
\begin{lemma}\label{le1} \cite{Majda 2001} (Product estimates)
For $m\in \mathbb{Z}^+$ and $m\geq2$, we have
\bal
&||\mathbf{f}\mathbf{g}||_{H^m}\leq C||\mathbf{f}||_{H^m}||\mathbf{g}||_{H^m},\label{cj1}
\\&||\mathbf{f}\mathbf{g}||_{H^m}\leq C\big(||\mathbf{f}||_{L^\infty}+||\na^m\mathbf{f}||_{L^\infty}\big)||\mathbf{g}||_{H^m}.\label{cj2}
\end{align}
\end{lemma}
\begin{lemma}\label{lem3.1} Under the assumptions of Theorem \ref{the1.1}, for all $M\in\mathbb{Z}^+\cup\{0\}$, the following estimates hold
\bal\label{l}
E(t)\triangleq||f,g,h||_{H^3}\leq Ce^{-t}\ep\big(||a_0||_{L^2}||\hat{a}_0||_{L^1}+||m_0||_{L^2}||\hat{m}_0||_{L^1}+||m_0||_{L^2}\big)
\end{align}
and
\bal\label{lj}
||U,V||_{W^{M,\infty}}\leq Ce^{-t}||\hat{a}_0,\hat{m}_0||_{L^1}.
\end{align}
\end{lemma}
{\bf Proof of Lemma \ref{lem3.1}}\quad  Recalling that $h=-(\pa_1+\pa_2)m$ and $m=e^{t\Delta}m_0$, then using the condition $\mathrm{supp}\ \hat{m}_0(\xi)\subset\mathcal{C}$, one has
\bal\label{h}
||h||_{H^M}=||(\pa_1+\pa_2)m||_{H^M}\leq Ce^{-t}\ep||m_0||_{L^2}.
\end{align}
Direct calculations show that
\bbal
f^1&=-U\cd\na U^1-V\cd\na V^1-V^1\div V
\\&=\pa_1a\pa_2\pa_2a-\pa_2a\pa_1\pa_2a-2m(\pa_1+\pa_2)m
\\&=(\pa_1+\pa_2)a\pa_2\pa_2a-\pa_2a\pa_2(\pa_1+\pa_2)a-2m(\pa_1+\pa_2)m,
\\f^2&=-U\cd\na U^2-V\cd\na V^2-V^2\div V
\\&=-\pa_1a\pa_2\pa_1a+\pa_2a\pa_1\pa_1a-2m(\pa_1+\pa_2)m
\\&=-(\pa_1+\pa_2)a\pa_1\pa_2a+\pa_2a\pa_1(\pa_1+\pa_2)a-2m(\pa_1+\pa_2)m,
\\g^1&=-U\cd\na V^1-V\cd\na U^1
\\&=\pa_1a\pa_2m-\pa_2a\pa_1m-m(\pa_1+\pa_2)\pa_2a
\\&=(\pa_1+\pa_2)a\pa_2m-\pa_2a(\pa_1+\pa_2)m-m(\pa_1+\pa_2)\pa_2a,
\\g^2&=-U\cd\na V^2-V\cd\na U^2
\\&=\pa_1a\pa_2m-\pa_2a\pa_1m+m(\pa_1+\pa_2)\pa_1a
\\&=(\pa_1+\pa_2)a\pa_2m-\pa_2a(\pa_1+\pa_2)m+m(\pa_1+\pa_2)\pa_1a.
\end{align*}
For all $M\in\mathbb{Z}^+\cup\{0\}$, using the fact $||\mathbf{f}||_{L^\infty}\leq C||\hat{\mathbf{f}}||_{L^{1}}$, we deduce
\bal\label{a}
||\nabla^M a||_{L^\infty}\leq Ce^{-t}\bi\||\xi|^M\hat{a}_0\bi\|_{L^{1}}\leq Ce^{-t}||\hat{a}_0||_{L^{1}}.
\end{align}
and
\bal\label{m}
||\nabla^M m||_{L^\infty}\leq Ce^{-t}\bi\||\xi|^Me^{-(|\xi|^2-1)t}\hat{m}_0\bi\|_{L^{1}}\leq Ce^{-t}||\hat{m}_0||_{L^{1}}.
\end{align}
where we have used the conditions $\mathrm{supp}\ \hat{a}_0(\xi)\subset \mathcal{C}$ and $\mathrm{supp}\ \hat{m}_0(\xi)\subset \mathcal{C}$.

By the product estimate \eqref{cj2}, then we have
\bal\label{f1}
||f^1||_{H^M}\leq&~ C\big(||(\pa_1+\pa_2)a||_{H^{M+1}}||a||_{W^{M+2,\infty}}+||(\pa_1+\pa_2)m||_{H^M}||m||_{W^{M,\infty}}\big)\nonumber\\
\leq&~ Ce^{-t}\ep\big(||a_0||_{L^2}||\hat{a}_0||_{L^1}+||m_0||_{L^2}||\hat{m}_0||_{L^1}\big).
\end{align}
An argument similar to that used above, we also have
\bal\label{f2}
||f^2,g^1,g^2||_{H^M}\leq&~ Ce^{-t}\ep\big(||a_0||_{L^2}||\hat{a}_0||_{L^1}+||m_0||_{L^2}||\hat{m}_0||_{L^1}\big).
\end{align}
Combining \eqref{h}, \eqref{f1} and \eqref{f2} yields the desired result \eqref{l}. \eqref{lj} is just a direct consequence of \eqref{a} and \eqref{m}. Thus, we end the proof of Lemma \ref{lem3.1}. $\Box$
\section{Proof of Theorem \ref{the1.1}}
By standard energy method, we can obtain that there exists a unique solution $(u,v,\theta)$ to \eqref{2d-tcm} on some time interval $[0,T^*)$, where $T^*$ is the maximal time of existence of solution $(u,v,\theta)$. It remains to prove $T^*=\infty$.

{\bf Step 1: The Estimate of $w$}.\\
Applying $D^\ell$ to the both sides of Eq.$\eqref{c-tcm}_1$, taking the inner product with $D^\ell w$ and summing the resulting over $|\alpha|\leq 3$, we get
\bal\label{es-wc}
\fr12\frac{\dd}{\dd t}||w||^2_{H^3}+||w||^2_{H^3}\triangleq\sum^{5}_{i=1}I_i,
\end{align}
where
\bbal
&I_1=-\sum_{0<|\ell|\leq 3}\int_{\R^2}[D^{\ell},w\cd] \na w\cd D^\ell w\dd x,\qquad I_2=-\sum_{0<|\ell|\leq 3}\int_{\R^2}[D^{\ell},U\cd] \na w\cd D^{\ell}w\dd x,\\
&I_3=-\sum_{0\leq |\ell|\leq 3}\int_{\R^2}D^{\ell}\big(c\cd\na c+c\div c\big)\cd D^\ell w\dd x,\qquad I_4=-\sum_{0\leq|\ell|\leq 3}\int_{\R^2}D^{\ell}f\cd D^\ell w\dd x\\
&I_5=-\sum_{0\leq|\ell|\leq 3}\int_{\R^2}D^{\ell}\big(c\cd\na V+c\div V+V\cd\na c+V\div c+w\cd\na U\big)\cd D^{\ell}w\dd x.
\end{align*}
Next, we need to estimate the above terms one by one.

For the terms $I_1$ and $I_2$, notice that the embedding $H^2(\R^2)\hookrightarrow L^\infty(\R^2)$, using the commutate estimate \eqref{jhz1} and \eqref{jhz2}, respectively, one has
\bal
I_1\leq&~\sum_{0<|\ell|\leq 3}||[D^{\ell},w\cd] \na w||_{L^2}||\na w||_{H^2}\nonumber\\
\leq&~C||\na w||_{L^\infty}||w||_{H^3}||\na w||_{H^2}\nonumber\\
\leq&~CA^\frac12(t)B(t),\label{z1}\\
I_2\leq&~\sum_{0<|\ell|\leq 3}||[D^{\ell},U\cd] \na w||_{L^2}||\na w||_{H^2}\nonumber\\
\leq& ~C\big(||\na U||_{L^\infty}+||\na^3U||_{L^\infty}\big)||\na w||^2_{H^2}\nonumber\\
\leq&~ C||\na U,\na^3U||_{L^\infty}A(t).\label{z3}
\end{align}
By the product estimate \eqref{cj1}, one has
\bal\label{z2}
I_3\leq&~ ||c\cd\na c+c\div c||_{H^3}||w||_{H^3}\nonumber\\
\leq&~ C||\na c||_{H^3}||c||_{H^3}||w||_{H^3}\nonumber\\
\leq&~CA^\frac12(t)B(t).
\end{align}
Using H\"{o}lder's inequality gives
\bal
I_4\leq&~||f||_{H^3}||w||_{H^3}\leq CE(t)+E(t)A(t).\label{z6}
\end{align}
For the last term $I_5$, by product estimate \eqref{cj2}, one has
\bal
I_5\leq&~||c\cd\na V+c\div V+V\cd\na c+V\div c+w\cd\na U||_{H^3}||w||_{H^3}\nonumber\\
\leq&~ C||\nabla U,\nabla V,\nabla^4 U,\nabla^4 V||_{L^\infty}||w,c||^2_{H^3}+C||V,\nabla^3 V||_{L^\infty}||\na c||_{H^3}||w||_{H^3}\nonumber\\
\leq&~ C||\nabla U,\nabla^4 U,V,\nabla V,\nabla^3 V,\nabla^4 V||_{L^\infty}A(t)+\fr12||\na c||^2_{H^3}.\label{z5}
\end{align}
Gathering the above estimates $I_1-I_5$ to \eqref{es-wc} together yields
\bal\label{wc}
\frac{\dd}{\dd t}||w||^2_{H^3}+||w||^2_{H^3}\leq CA^\frac12(t)B(t)+C\big(||U,V||_{W^{4,\infty}}+E(t)\big)A(t)+CE(t),
\end{align}
{\bf Step 2: The Estimates of $c$ and $\theta$}.\\
Applying $D^\ell$ to the both sides of Eqs.$\eqref{c-tcm}_2$ and $\eqref{c-tcm}_3$, taking the inner product with $D^\ell c$ and $D^\ell \theta$, respectively, then summing the resulting over $|\alpha|\leq 3$, we get
\bal\label{es-cc}
\fr12\frac{\dd}{\dd t}\Big(||c||^2_{H^3}+||\theta||^2_{H^3}\Big)+||\na c||^2_{H^3}\triangleq\sum^{6}_{i=1}J_i,
\end{align}
where
\bbal
&J_1=-\sum_{0<|\ell|\leq 3}\int_{\R^2}[D^{\ell},w\cd] \na c\cd D^\ell c\dd x-\sum_{0<|\ell|\leq 3}\int_{\R^2}[D^{\ell},w\cd] \na \theta\cd D^\ell \theta\dd x,
\\
&J_2=-\sum_{0<|\ell|\leq 3}\int_{\R^2}[D^{\ell},U\cd] \na \theta\cd D^{\ell}\theta\dd x-\sum_{0<|\ell|\leq 3}\int_{\R^2}[D^{\ell},U\cd] \na c\cd D^{\ell}c\dd x,
\\
&J_3=-\sum_{0\leq |\ell|\leq 3}\int_{\R^2}D^{\ell}(c\cd\na w)\cd D^\ell c\dd x,\quad J_4=-\sum_{0\leq |\ell|\leq 3}\int_{\R^2}D^{\ell}(V\cd\na w)\cd D^\ell c\dd x,\\
&J_5=-\sum_{0\leq |\ell|\leq 3}\int_{\R^2}D^{\ell}[c\cd\na U+w\cd\na V]\cd D^\ell c\dd x,\\
&J_6=-\sum_{0\leq|\ell|\leq 3}\int_{\R^2}D^{\ell}g\cd D^\ell c\dd x-\sum_{0\leq|\ell|\leq 3}\int_{\R^2}D^{\ell}h\cd D^\ell \theta\dd x.
\end{align*}
Next, we need to estimate the above terms one by one.

Using the commutate estimate \eqref{jhz1} and \eqref{jhz2}, respectively, we obtain
\bal
J_1\leq&~\sum_{0<|\ell|\leq 3}|||[D^{\ell},w\cd]\na c||_{L^2}||\na c||_{H^2}+\sum_{0<|\ell|\leq 3}|||[D^{\ell},w\cd]\na \theta||_{L^2}||\na \theta||_{H^2}\nonumber\\
\leq&~C\big(||\na c||_{H^2}||\na w||_{L^\infty}+||w||_{H^3}||\na c||_{L^\infty}\big)||\na c||_{H^2}\nonumber\\
\quad&+C\big(||\na \theta||_{H^2}||\na w||_{L^\infty}+||w||_{H^3}||\na \theta||_{L^\infty}\big)||\na \theta||_{H^2}\nonumber\\
\leq&~CA^\frac12(t)B(t),\label{j1}\\
J_2\leq&~\sum_{0<|\ell|\leq 3}||[D^{\ell},U\cd] \na c||_{L^2}||\na c||_{H^2}+\sum_{0<|\ell|\leq 3}||[D^{\ell},U\cd] \na \theta||_{L^2}||\na \theta||_{H^2}\nonumber\\
\leq& ~C\big(||\na U||_{L^\infty}+||\na^3U||_{L^\infty}\big)||\na w,\na \theta||^2_{H^2}\nonumber\\
\leq&~ C||\na U,\na^3U||_{L^\infty}A(t),\label{j2}
\end{align}
Due to the fact $c\cd\na w=\div(w\otimes c)-w\div c$, integrating by parts, one has
\bal
J_3=&~\sum_{0\leq |\ell|\leq 3}\int_{\R^2}D^{\ell}(w\otimes c)\cd D^\ell \nabla c\dd x+\sum_{0\leq |\ell|\leq 3}\int_{\R^2}D^{\ell}(w\div c)\cd D^\ell c\dd x\nonumber\\
\leq&~ C||w||_{H^3}||c||_{H^3}||\na c||_{H^3}
\nonumber\\
\leq&~ CA(t)B(t)+\fr14||\na c||^2_{H^3}.\label{z4}
\end{align}
Similarly, we deduce
\bal
J_4=&~\sum_{0\leq |\ell|\leq 3}\int_{\R^2}D^{\ell}(w\otimes V)\cd D^\ell \nabla c\dd x+\sum_{0\leq |\ell|\leq 3}\int_{\R^2}D^{\ell}(w\div V)\cd D^\ell c\dd x\nonumber\\
\leq&~ C||w\otimes V||_{H^3}||\na c||_{H^3}+C||w\div V||_{H^3}||c||_{H^3}\nonumber\\
\leq&~ C||V,\na^3V||_{L^\infty}||c||_{H^3}||\na c||_{H^3}+C||\nabla V,\na^4V||_{L^\infty}||w||_{H^3}||c||_{H^3}\nonumber\\
\leq&~ C\big(||\nabla V,\na^4V||_{L^\infty}+||V,\na^3V||^2_{L^\infty}\big)A(t)+\fr14||\na c||^2_{H^3}.
\end{align}
For the term $J_5$, by the product estimate \eqref{cj2}, one has
\bal
J_5\leq&~||c\cd\na U+w\cd\na V||_{H^3}||c||_{H^3}\leq~ C||\nabla U,\nabla V,\nabla^4 U,\nabla^4 V||_{L^\infty}A(t).
\end{align}
Using H\"{o}lder's inequality gives
\bal
J_6\leq&~||g||_{H^3}||c||_{H^3}+||h||_{H^3}||\theta||_{H^3}\leq CE(t)+E(t)A(t).
\end{align}
Gathering the above estimates $J_1-J_6$ to \eqref{es-cc} together yields
\bal\label{cc}
\frac{\dd}{\dd t}||c,\theta||^2_{H^3}+||\na c||^2_{H^3}\leq& C\big(A^\frac12(t)+A(t)\big)B(t)+CE(t)\nonumber\\
&~+C\big(||U,V||_{W^{4,\infty}}+||V,\na^3V||^2_{L^\infty}+E(t)\big)A(t).
\end{align}
{\bf Step 3: The Estimate of Crossing Term $\sum_{0\leq |\ell|\leq 2}\langle D^\ell c,D^\ell\na \theta\rangle$}.\\
Applying $D^\ell$ to $\eqref{c-tcm}_2$ and $\eqref{c-tcm}_3$, taking the scalar product of them with $D^\ell \na \theta$ and $D^\ell \div c$, respectively, adding them together and then summing the resulting over $|\alpha|\leq 2$, we get
\bal\label{es-theta}
\frac{\dd}{\dd t}\sum_{0\leq |\ell|\leq 2}\langle D^\ell c,D^\ell\na \theta\rangle+||\na \theta||^2_{H^2}-||\div c||^2_{H^2}-\sum_{0\leq |\ell|\leq 2}\langle D^\ell \De c,D^\ell\na \theta\rangle\triangleq\sum^{3}_{i=1}K_i,
\end{align}
where
\bbal
&K_1=-\sum_{0\leq |\ell|\leq 2}\int_{\R^2}D^{\ell}[w\cd\na c+c\cd\na w]\cd D^\ell \na\theta\dd x+\sum_{0\leq |\ell|\leq 2}\int_{\R^2}D^{\ell}(w\cd\na \theta)\cd D^\ell \div c\dd x,
\\&K_2=-\sum_{0\leq |\ell|\leq 2}\int_{\R^2}D^\ell[V\cd\na w+c\cd\na U+U\cd\na c+w\cd\na V]\cd D^\ell\na \theta\dd x
\\&\qquad+\sum_{0\leq |\ell|\leq 2}\int_{\R^2}D^{\ell}(U\cd\na \theta)\cd D^\ell \div c\dd x,
\\&K_3=\sum_{0\leq |\ell|\leq 2}\int_{\R^2}D^{\ell}g\cd D^\ell \na\theta\dd x-\sum_{0\leq |\ell|\leq 2}\int_{\R^2}D^{\ell}h\cd D^\ell \div c\dd x.
\end{align*}
By Leibniz's formula and H\"{o}lder's inequality, we have
\bal
K_1\leq&~||w\cd \na c+c\cd\na w||_{H^2}||\na\theta||_{H^2}+||w\cd \na \theta||_{H^2}||\na c||_{H^2}\nonumber\\
\leq&~ C||w||_{H^3}||c||_{H^3}||\na\theta||_{H^2}\nonumber\\
\leq&~CA(t)B(t)+\fr12 ||\na\theta||_{H^2},\label{z6}\\
K_2\leq&~||V\cd\na w+c\cd\na U+U\cd\na c+w\cd\na V||_{H^2}||\na\theta||_{H^2}+||U\cd \na \theta||_{H^2}||\na c||_{H^2}\nonumber\\
\leq&~ C||U,V||_{W^{3,\infty}}||c||_{H^3}||\theta||_{H^3}\nonumber\\
\leq&~C||U,V||_{W^{3,\infty}}A(t),\label{z7}
\\K_3\leq&~C||g,h||_{H^2}||c,\theta||_{H^3}
\leq CE(t)+E(t)A(t).\label{z8}
\end{align}
Gathering the above estimates $K_1-K_6$ to \eqref{es-theta} together yields
\bal\label{theta}
\frac{\dd}{\dd t}\sum_{0\leq |\ell|\leq 2}\langle D^\ell c,D^\ell\na \theta\rangle&+\fr12||\na \theta||^2_{H^2}-||\div c||^2_{H^2}-\sum_{0\leq |\ell|\leq 2}\langle D^\ell \De c,D^\ell\na \theta\rangle\nonumber\\
&~\leq CA(t)B(t)+C\big(||U,V||_{W^{3,\infty}}+E(t)\big)A(t)+CE(t)
\end{align}
{\bf Step 4: Closure of The A Priori Estimates}.\\
By simple computations, we deduce easily that for some suitable positive constant $\gamma$
$$
A(t)+\gamma \sum_{0\leq |\ell|\leq 2}\langle D^\ell v,D^\ell\na \theta\rangle\thickapprox A(t)
$$
and
$$||w||^2_{H^3}+||\na c||^2_{H^3}+\frac{\gamma}{2}||\na \theta||^2_{H^2}-\gamma||\div c||^2_{H^2}-\gamma\sum_{0\leq |\ell|\leq 2}\langle D^\ell \De c,D^\ell\na \theta\rangle\thickapprox B(t)
$$
Multiplying the inequality \eqref{es-theta} by $\gamma$, combining \eqref{wc} and \eqref{cc}, then integrating in time yields
\bal\label{es-con}
A(t)+\int_0^tB(s)\dd s\leq&~C\int_0^t\big(A^{\frac12}(s)+A(s)\big)B(s)\dd s+\int_0^tE(s)\dd s\nonumber\\
&+ C\int_0^t\big(||U,V||_{W^{4,\infty}}+||V,\nabla^3V||^2_{L^\infty}+E(s)\big)A(s)\dd s.
\end{align}
Let us define
\bbal
\Gamma\triangleq\sup\bi\{t\in[0,T^*): \sup_{\tau\in[0,t]}A(\tau)\leq \eta\ll1\bi\},
\end{align*}
where $\eta$ is a small enough positive constant which will be determined later on.

Assume that $\Gamma<T^*$. Choosing $\eta$ small enough such that the first term of RHS of \eqref{es-con} is absorbed, then we infer from \eqref{es-con} that for all $t\in[0,\Gamma]$
\bal\label{es-con1}
A(t)+\int_0^tB(s)\dd s\leq&~C\int_0^tE(s)\dd s+ C\int_0^t\big(||U,V||_{W^{4,\infty}}+||V,\nabla^3V||^2_{L^\infty}+E(s)\big)A(s)\dd s.
\end{align}
Then by Gronwall's inequality and \eqref{condition}, \eqref{l}--\eqref{lj}, we have for all $t\in[0,\Gamma]$
\bal\label{es-con2}
A(t)\leq&~C\bi(A_0+\int_0^tE(s)\dd s\bi) \exp\bi\{C\int_0^t\big(||U,V||_{W^{4,\infty}}+||V,\nabla^3V||^2_{L^\infty}+E(s)\big)\dd s\bi\}\nonumber\\
\leq&~C\varepsilon_0.
\end{align}
Choosing $\eta=2C\varepsilon_0$, thus we can get
\bbal
\sup_{\tau\in[0,t]}A(\tau)&\leq \fr\eta2 \quad\mbox{for}\quad t\leq \Gamma.
\end{align*}

So if $\Gamma<T^*$, due to the continuity of the solutions, we can obtain that there exists $0<\epsilon\ll1$ such that
\bbal
\sup_{\tau\in[0,t]}A(\tau)&\leq \eta \quad\mbox{for}\quad t\leq \Gamma+\epsilon<T^*,
\end{align*}
which is contradiction with the definition of $\Gamma$.

Thus, we can conclude $\Gamma=T^*$ and
\bbal
\sup_{\tau\in[0,t]}A(\tau)&\leq C<\infty \quad\mbox{for all}\quad t\in(0,T^*),
\end{align*}
which implies that $T^*=+\infty$. This completes the proof of Theorem \ref{the1.1}. $\Box$

\section*{Acknowledgments} J. Li is supported by NSFC (No.11801090). Y. Yu is supported by NSF of Anhui Province (No.1908085QA05).

\end{document}